\documentstyle[12pt]{article}

\makeatletter
\@ifundefined{url}{\def\url#1{\texttt{#1}}}{}
\makeatother

\begin{document}
\newfont{\blb}{msbm10 scaled\magstep1} 

\newtheorem{theo}{Theorem}[section]

\newtheorem{defi}[theo]{Definition}
\newtheorem{prop}[theo]{Proposition}
\newtheorem{lemm}[theo]{Lemma}
\newtheorem{coro}[theo]{Corollary}
\pagestyle{myheadings}
\date{}
\author{Layla Sorkatti \\
Gunnar Traustason \\
Department of Mathematical Sciences \\ 
University of Bath, UK}
\title{Nilpotent symplectic alternating algebras I}
\maketitle
\begin{abstract}
\mbox{}\\
We develop a structure theory for nilpotent symplectic alternating algebras.\\\\
{\small Keywords: Non-associative algebras, Symplectic, Engel, Nilpotent, Powerful, p-group.}
\end{abstract}
\section{Introduction}
Symplectic alternating algebras have arisen from the study of $2$-Engel groups
(see~\cite{moravec,gt-2008}) but seem also to be of interest in their own right, with many
beautiful properties. Some general theory was developed in~\cite{saa,to,sor-thesis}. \\ 

\noindent
{\bf Definition}.
Let $F$ be a field. A {\it symplectic alternating algebra} over $F$ is a
triple $(L,(\ ,\ ),\cdot)$ where $L$ is a symplectic vector space
over $F$ with respect to a non-degenerate alternating form $(\ ,\ )$
and $\cdot$ is a bilinear and alternating binary operation on $L$
such that
      $$(u\cdot v,w)=(v\cdot w,u)$$
for all $u,v,w\in L$. \\ \\
Notice that $(u\cdot x,v)=(x\cdot v,u)=-(v\cdot x,u)=(u,v\cdot x)$ and thus
the multiplication from the right by $x$ is self-adjoint with respect
to the alternating form. As the alternating form is non-degenerate, $L$ is
of even dimension and we can pick a basis $x_{1},y_{1},\ldots ,x_{n},y_{n}$
with the property that $(x_{i},x_{j})=(y_{i},y_{j})=0$ and $(x_{i},y_{j})=
\delta_{ij}$ for $1\leq i\leq j\leq n$. We refer to a basis of this type as
a {\it standard basis}. \\ \\
Suppose we have any basis $u_{1},\ldots ,u_{2n}$ for $L$. The structure of $L$
is then determined from 
       $$(u_{i}u_{j},u_{k})=\gamma_{ijk},\ \ 1\leq i<j<k\leq 2n.$$
The map $L^{3}\rightarrow F,\,(u,v,w)\mapsto (u\cdot v,w)$ is an alternating 
ternary form and each alternating ternary form on a given symplectic vector 
space, with a non-degenerate alternating form, defines a unique symplectic alternating algebra. Classifying symplectic
alternating algebras of dimension $2n$ over a field $F$ is then equivalent to
finding all the $\mbox{Sp}\,(V)$-orbits of $(\wedge^{3}V)^{*}$ under the natural
action, where $V$ is a symplectic vector space of dimension $2n$ with a 
non-degenerate alternating form. 
Suppose that $F$ is a finite field and suppose that the disjoint 
$\mbox{Sp\,}(V)$-orbits of $(\bigwedge^{3}V)^{*}$ are 
$u_{1}^{\mbox{Sp}(V)}$, $\ldots$, $u_{m}^{\mbox{Sp}(V)}$. Then 
    $$m\leq |F|^{2n\choose 3}=|(\wedge^{3}\,V)^{*}|\leq m |\mbox{Sp}(V)|
      \leq m|F|^{2n+1\choose 2}.$$
It follows that $m=|F|^{\frac{4n^{3}}{3}+O(n^{2})}$. Because of the sheer growth,
a general classification of symplectic alternating algebras seems impossible.
There is a close connection
between symplectic alternating algebras over the field $\mbox{GF}(3)$ of three elements and
a certain class of $2$-Engel groups and in~\cite{saa} the symplectic alternating
algebras over $\mbox{GF}(3)$ of dimension up to $6$ were classified. There
are $31$ such algebras of dimension $6$ of which $15$ are simple. 
We would like to mention here also the work of Atkinson~\cite{Ank} who
in his thesis looked at alternating ternary forms over $\mbox{GF\,}(3)$ in 
order to study  a certain class of groups of exponent $3$. \\ \\
As we said before, some general theory was developed in~\cite{saa,to}. In 
particular a well-known dichotomy property for Lie algebras also
holds for symplectic alternating algebras. Thus a symplectic alternating algebra
is either semi-simple or has a non-trivial abelian ideal. Another interesting 
property is that any symplectic
alternating algebra that is nilpotent-by-abelian must be nilpotent. \\ \\ 
In this paper and its sequel we study nilpotent symplectic alternating algebras. This paper deals with the structure theory and also gives the classification of nilpotent symplectic alternating algebras of dimension up to $8$ over any field. The sequel will be mostly about extending the classification of algebras to 
dimension $10$ that is far more involved than the algebras of lower
dimension.
Additionally, we introduce a special class of powerful p-groups that we call powerfully nilpotent groups.
These groups represent finite p-groups characterized by possessing a central series of a unique nature. Turning back to this paper, we will in Section $2$ describe some general results that in 
particular
lead to specific type of presentations that we call nilpotent presentations.
All algebras with a nilpotent presentation are nilpotent and conversely any
nilpotent algebra will have a nilpotent presentation. In Section $3$ we will 
focus
on the algebras that are of maximal class and we will see that their structure
is very rigid. Finally to illustrate the theory we will in Section $4$ classify
the nilpotent symplectic alternating algebras of dimension up to $8$ over
an arbitrary field $F$. \\ \\
In this paper we will adopt the {\it left-normed} convention for products. Thus $u_{1}u_{2}
\ldots u_{n}$ stands for $(\ldots (u_{1}u_{2})\cdots )u_{n}$. Also $U\leq V$
stands for `$U$ is a subspace of $V$'. \\ \\
Many of the terms that we use in this paper are analogous to the corresponding
terms for related structures. Thus a subspace $I$ of a symplectic alternating
algebra $L$ is an ideal if $IL\leq I$. From [5] we know that $I^{\perp}$ is an
ideal whenever $I$ is an ideal.
The definition of a nilpotent symplectic alternating algebra causes no problem
either. \\ \\
{\bf Definition} A symplectic alternating algebra $L$ is {\it nilpotent} if there 
exists an ascending
chain of ideals $I_{0},\ldots, I_{n}$ such that 
    $$\{0\}=I_{0}\leq I_{1}\leq \cdots \leq I_{n}=L$$
and $I_{s}L\leq I_{s-1}$ for $s=1,\ldots ,n$. 
The smallest possible $n$ is then called the {\it nilpotence class} of $L$. \\ \\
{\bf Definition}. More generally, if $I_{0}\leq I_{1}\leq \ldots \leq I_{n}$ is 
any chain of ideals of $L$ then we say that this chain is central in $L$ if
$I_{s}L\leq I_{s-1}$ for $s=1,\ldots ,n$. \\ \\
We define the lower and upper central series in an analogous way to
related structures like associative algebras and Lie algebras. Thus we define
the lower central series recursively by $L^{1}=L$ and $L^{n+1}=L^{n}L$, and
the upper central series by $Z_{0}(L)=\{0\}$ and $Z_{n+1}(L)=
\{x\in L:\,xL\subseteq Z_{n}(L)\}$. It is readily seen that the terms of the
lower and the upper central series are all ideals of $L$. The following beautiful property was proved
in [5] and will be used frequently
\begin{equation}
Z_{n}(L)=(L^{n+1})^{\perp}.
\end{equation}
\mbox{}\\
{\bf Remark}. Notice however that the lack of the Jacobi identity means that 
many properties that hold for Lie algebras do not hold for symplectic 
alternating algebras. As the following example shows, it is not true in general 
that the product of two ideals is an ideal. That example also shows that
the formula $L^{i}L^{j}\leq L^{i+j}$ does not hold in general. \\ \\
{\bf Example}. Consider the $12$-dimensional symplectic alternating algebra
which has a standard basis $x_{1},y_{1},x_{2},y_{2},x_{3},y_{3},x_{4},y_{4},x_{5},y_{5},x_{6},y_{6}$ where
     $$(x_{3}y_{5},y_{6})=(x_{2}y_{4},y_{6})=(x_{1}y_{4},y_{5})=
  (y_{1}y_{2},y_{3})=1$$
and $(uv,w)=0$ if $u,v,w$ are basis elements where $\{u,v,w\}\not\in 
\{\{x_{3},y_{5},y_{6}\}$ ,
$\{x_{2},y_{4},y_{6}\},\{x_{1},y_{4},y_{5}\},\{y_{1},y_{2},y_{3}\}\}$. Notice
that this implies that  
$$\begin{array}{lll}
    x_{3}y_{5}=x_{6}, & x_{1}y_{4}=x_{5}, & y_{2}y_{3}=x_{1}, \\
    x_{3}y_{6}=-x_{5}, & x_{1}y_{5}=-x_{4}, & y_{4}y_{5}=-y_{1}, \\
    x_{2}y_{4}=x_{6}, & y_{1}y_{2}=x_{3}, & y_{4}y_{6}=-y_{2}, \\
  x_{2}y_{6}=-x_{4}, & y_{1}y_{3}=-x_{2}, & y_{5}y_{6}=-y_{3}. 
\end{array}$$
From this one sees that 
\begin{eqnarray*}
   L^{2} & = & Fx_{6}+Fx_{5}+\cdots
        +Fx_{1}+Fy_{1}+Fy_{2}+Fy_{3}, \\
  L^{3} & = & Fx_{6}+Fx_{5}+\cdots +
Fx_{1}, \\
  L^{4} & = & Fx_{6}+Fx_{5}+Fx_{4}, \\
 L^{5} & = & 0, \\
 L^{2}L^{2} & = & Fx_{3}+Fx_{2}+Fx_{1}.
\end{eqnarray*}
In particular $L$ is nilpotent of class $4$, $L^{2}L^{2}$ is not an ideal
and $L^{2}L^{2}\not\leq L^{4}$. \\ \\
This example indicates that symplectic alternating algebras do differ
from Lie algebras. We are going to see in the following sections that there
are some shared properties but the next lemma underlines the difference by 
showing that the two classes of algebras do not have many algebras in common 
when the characteristic is not $2$. In fact only the symplectic alternating
algebras that are obviously Lie algebras are there, namely those
of class  at most $2$. 
\begin{lemm} Let $L$ be a symplectic alternating algebra where $\mbox{char\,}L
\not =2$ and  $L$ is either associative or a Lie algebra. Then $L^{3}=\{0\}$. 
\end{lemm}
{\bf Proof}. Let us first assume that $L$ is associative. We then have
   $$0=(xyz-x(yz),t)=(x,tzy-t(yz))=(x,tzy-tyz)$$
for all $x,y,z,t\in L$. It follows that $tzy=tyz=-ytz$ for all $t,z,y\in L$.
Using this last property repeatedly we get that 
     $$xyz=-zxy=yzx=-xyz$$
and thus $2xyz=0$ for all $x,y,z\in L$. As \mbox{char\,}$L\not =2$, it 
follows that $L^{3}=0$. \\ \\
Now suppose $L$ is a Lie algebra. We then have
   $$0=(xyz+yzx+zxy,t)=(x,tzy-t(yz)-tyz)=2(x,tzy-tyz).$$
As \mbox{char\,}$L\not =2$, it follows again that $tzy=tyz$
for all $t,z,y\in L$ and this implies again that $L^{3}=\{0\}$. $\Box$. \\ \\
One handicap that the symplectic alternating algebras have is that when $I$ is
an ideal then $L/I$ is in general only an alternating algebra as there is no
natural way of inducing an alternating form on this quotient. For example simply for the reason that the quotient can have odd dimension. There is however a
weaker form of a quotient structure that we can associate to any ideal
$I$ of $L$ that works. Thus for any ideal $I$ we have that $(I^{\perp}+I)/I$
is a well defined symplectic alternating algebra with the natural induced multiplciaton and where the induced alternating form is given by $(u+I,v+I)=(u,v)$
for $u,v\in I^{\perp}$. The reader can easily convince himself that this is
well defined and that $((I^{\perp}+I)/I)^{\perp}=0$. This algebra is also 
isomorphic to $I^{\perp}/(I\cap I^{\perp})$ that has a similar naturally induced
structure as a symplectic alternating algebra. \\ \\
{\bf Remark}. There are some familiar facts for Lie algebras that do not
rely on the Jacobi identity and remain true for symplectic alternating algebras.
Such properties are particularly useful as we can use them when dealing
with quotients $L/I$ where we only know that the resulting algebra is 
alternating. For example $L^{2}$ has co-dimension at least $2$ in any nilpotent alternating algebra
$L$ of dimension greater than or equal to $2$. From this and the duality given in (1), it
follows immediately that the dimension of $Z(L)$ is at least $2$ for any
non-trivial nilpotent symplectic algebra which is something that we will also see later
as a corollary of Lemma 2.1. 

\section{General Structure Theory}
 We next see
that, like for Lie algebras, all minimal sets of generators have the
same number of elements and we can thus introduce the notion of a rank.  \\ \\
{\bf Definition}. Let $L$ be a nilpotent symplectic alternating algebra. 
We say that $\{x_{1},\ldots ,x_{r}\}$ is a {\it minimal set of generators} if
these generate $L$ (as an algebra) and no proper subset generates $L$. \\
\begin{lemm} Let $L$ be a nilpotent symplectic alternating algebra. Any minimal
set of generators has the same size which is $\mbox{dim\,}L-\mbox{dim\,}L^{2}$.
\end{lemm}
{\bf Proof}\ \ Let $x_{1},\ldots ,x_{r}\in L$ and let $M$ be the subalgebra
of $L$ generated by these elements. It suffices to show that $L=M$ if and
only if $x_{1}+L^{2}, \cdots ,x_{r}+L^{2}$ generate $L/L^{2}$ as a vector
space. Suppose first that $L=M$. Notice that $M=Fx_{1}+\cdots +Fx_{r}+
M\cap L^{2}$ and thus it is clear that $L/L^{2}$ is generated by
$x_{1}+L^{2},\ldots ,x_{r}+L^{2}$ as a vector space. Conversely suppose
now that the images of $x_{1},\ldots ,x_{r}$ in  $L/L^{2}$ generate $L/L^{2}$ as
a vector space. An easy induction shows that 
      $$L=M+L^{s+1},\ \ L^{s}=M^{s}+L^{s+1}$$
for all integers $s\geq 1$. If the class of $L$ is $n$, we get in particular
that $L=M+L^{n+1}=M$. $\Box$ \\ \\
{\bf Definition}. Let $L$ be a nilpotent symplectic alternating
algebra. The unique smallest number of generators for $L$, as an algebra,
is called the {\it rank} of $L$ and is denoted $r(L)$. \\ \\
By last lemma we know that $r(L)=\mbox{dim\,}L-\mbox{dim\,}L^{2}$. This has
the following curious consequence.
\begin{coro} Let $L$ be a nilpotent symplectic alternating algebra.
We have $r(L)=\mbox{dim\,}Z(L)$. In particular if $L\not =\{0\}$ then 
$\mbox{dim\,}Z(L)\geq 2$.
\end{coro}
{\bf Proof}. From (1) we know that $Z(L)=(L^{2})^{\perp}$. Therefore 
  $$r(L)=\mbox{dim\,}L-\mbox{dim\,}L^{2}=
\mbox{dim\,}(L^{2})^{\perp}=\mbox{dim\,}Z(L).$$
Finally, we cannot have $r(L)=1$ as then we would have that $L$ is 
one-dimensional. Hence $\mbox{dim\,}Z(L)\geq 2$. $\Box$. 
\begin{lemm} 
Let $I,J$ be ideals of a nilpotent symplectic alternating algebra where $I\subseteq J$. If 
$\mbox{dim\,}J=\mbox{dim\,}I+1$ then $I\leq J$ is central. If $I$ is
an ideal such that $\mbox{dim\,}I<2n=\mbox{dim\,}L$ then there exists an ideal $J$
such that $\mbox{dim\,}J=\mbox{dim\,}I+1$. If furthermore $I$ is an isotropic
ideal and $\mbox{dim\,}I<n$ then $J$ can be chosen to be isotropic. 
\end{lemm}
{\bf Proof}\ \ Suppose $J=I+Fx$ for some $x\in L$. Let $y\in L$. To show
that $I\leq J$ is central, it suffices to show that $x\cdot y\in I$. Suppose
that $xy=u_{1}+ax$ for $u_{1}\in I$ and $a\in F$. As $I\unlhd L$ it follows by
induction that $xy^{r}=u_{r}+a^{r}x$ for some $u_{r}\in I$. If $L$ is nilpotent
of class at most $m$ it follows that $0=u_{m}+a^{m}x$ and hence $a=0$. \\ \\
For the latter part suppose first that $I$ is any ideal such that 
$\mbox{dim\,}I<2n$. Let $m$ be the largest positive integer such that 
$L^{m}\not \leq I$. Pick $u\in L^{m}\setminus I$. Then $J=I+Fu$ is the required ideal such that $I\leq J$ is central. Now suppose furthermore that $I$ is 
isotropic and that $\mbox{dim\,}I<n$. Then $I^{\perp}$ is also an ideal
and $I< I^{\perp}$. Let $m$ be the largest non-negative integer such that 
$I^{\perp}\underbrace{L\cdots L}_{m}\not\leq I$. Let $u\in I^{\perp}\underbrace{L\cdots L}_{m} \setminus I$ and again
the ideal $J=I+Fu$ is the one required. $\Box$ \\ \\
{\bf Remark}. Let $U,V$ be subspaces of $L$. Notice that 
    $$UV=0\Leftrightarrow (UV,L)=0 \Leftrightarrow
     (UL,V)=0\Leftrightarrow  UL\leq V^{\perp}.$$
In other words we have that $U$ annihilates $V$ if and only if it
annihilates $L/V^{\perp}$. This is a useful property that we will be making use of later. 
\begin{lemm} Let $I$ and $J$ be ideals in a symplectic alternating algebra
$L$. We have that $I\cdot L\leq J$ if and only if $I\cdot J^{\perp}=\{0\}$. In 
particular $L^{m}\cdot Z_{m}(L)=\{0\}$ for all $m\geq 1$. 
\end{lemm}
{\bf Proof}\ \ From the property given in last remark, we know that
$I$ annihilates $L/J$ if and only if $I$ annihilates $J^{\perp}$. The second
part follows from this, the fact that $L^{m}$ annihilates $L/L^{m+1}$, and 
the fact that $(L^{m+1})^{\perp}=Z_{m}(L)$. $\Box$ \\ \\
{\bf Remark}. It follows in particular that $I\cdot I^{\perp}=\{0\}$ for
any ideal $I$. In particular any isotropic ideal is  abelian. Notice also that the property $L^{m}Z_{m}(L)=\{0\}$ is equivalent to the fact that
$Z_{m}(L)$ annihilates $L/Z_{m-1}(L)$. 
\\ \\
{\bf Remark}. We have seen in the introduction that it is not true 
in general that $L^{i}L^{j}\leq L^{i+j}$. As $(L^{m})^{\perp}=Z_{m-1}(L)$ for 
all $m\geq 1$, we
however have that 
$$\begin{array}{l}
  L^{i}L^{j}\leq L^{i+j} \Leftrightarrow (L^{i}L^{j},Z_{i+j-1}(L))=0
\Leftrightarrow (L^{i}Z_{i+j-1}(L),L^{j})=0 \\
\mbox{\ \ \ \ \ \ \ \ \ \ \ \ \ \ \ \  }\Leftrightarrow L^{i}Z_{i+j-1}(L)\leq Z_{j-1}(L).
\end{array}$$
The obvious fact that $L^{m}L\leq L^{m+1}$ thus gives us the interesting
fact from last lemma that $L^{m}Z_{m}(L)=\{0\}$. \\
\begin{lemm} Let $I$ be an ideal of $L$. Then $IL\leq I^{\perp}$ if and
only if $I$ is abelian. 
\end{lemm}
{\bf Proof}\ \ We have that $I$ annihilates $I$ if and only if $I$ annihilates
$L/I^{\perp}$. $\Box$ \\ \\
{\bf Remark}. As $I$ is an ideal we have in fact that $IL\leq I^{\perp}$ if
and only if $IL\leq I\cap I^{\perp}$. Here $I\cap I^{\perp}$ is the `isotropic
part' of $I$. 
\begin{lemm} Let $L$ be a nilpotent symplectic alternating algebra with ideals
$I,J$ where $J=I+Fx+Fy$, $(x,y)=1$ and $Fx+Fy\leq I^{\perp}$. Then $JL\leq I$.
Furthermore if $I$ is isotropic then $J$ is abelian. 
\end{lemm}
{\bf Proof}\ \  As $J$ is an ideal of $L$ and as $(xt,x)=0$ for all $t\in L$ we have that $I+Fx$ is an
ideal of $L$. By Lemma 2.3 we have that $I\leq I+Fx$ is central. 
Similarly $I\leq I+Fy$ is central and thus $JL\leq I$. For the second part
notice that if $I$ is isotropic then $I=J\cap J^{\perp}$ thus $JL\leq 
I=J\cap J^{\perp}$ and by Lemma 2.5 it follows that $J$ is abelian. $\Box$
\begin{lemm} Let $L$ be a nilpotent symplectic alternating algebra. Every ideal
$I$ of dimension $2$ is contained in $Z(L)$. Equivalently, every ideal
of co-dimension $2$ must contain $L^{2}$. 
\end{lemm}
{\bf Proof}. The second statement is a trivial fact that holds in
all nilpotent alternating algebras. The first statement  is a consequence
of this and of the duality given by $I\leq Z(L)\Leftrightarrow L^{2}=Z(L)^{\perp}
\leq I^{\perp}$. $\Box$
\begin{lemm} Let $I,J$ be ideals of a symplectic alternating algebra $L$ and let $x\in L$. We have $Jx\leq I$ if and only if
$I^{\perp}x\leq J^{\perp}$. 
\end{lemm}
{\bf Proof}. We have that $Jx\leq I$ is equivalent to $(ux,v)=0$ for all $u\in J$ and all $v\in I^{\perp}$. But this is equivalent
to saying that $(vx,u)=0$ for all $v\in I^{\perp}$ and $u\in J=(J^{\perp})^{\perp}$ and this is the same as saying that
$I^{\perp}x\leq J^{\perp}$. $\Box$ 
\begin{prop} Let $L$ be a symplectic alternating algebra. No term of the
upper central series has co-dimension $1$. Equivalently, no term of 
the lower central series has dimension $1$.
\end{prop}
{\bf Proof}\ \ The first fact is a well-known fact about alternating algebras
and follows from the fact that if $A$ is an alternating algebra then $A/Z(A)$
cannot be one-dimensional. Now the interesting second statement is a consequence
of this and the duality $(L^{r})^{\perp}=Z_{r-1}(L)$. $\Box$ 
\mbox{}\\ \\
In particular we have that 
                 $$\{0\}=I_{0}\leq I_{1}\leq \cdots \leq I_{m}=L$$
is an ascending central chain if and only if
                 $$L=I_{0}^{\perp}\geq I_{1}^{\perp}\geq \cdots \geq I_{m}^{\perp}=\{0\}$$
is a descending central chain. \\ \\
{\bf Remark}. Suppose that $L$ is any nilpotent alternating algebra such
that $L/L^{2}$ is $2$-dimensional. Then it follows immediately that 
the dimension of $L^{2}/L^{3}$ is at most $1$ and that the dimension
of $L^{3}/L^{4}$ is at most $2$. Using this general fact and Proposition 2.9
one can quickly show that all nilpotent symplectic alternating algebras of
dimension up to $4$ must be abelian. This is clear when the 
dimension is $2$. Now suppose that $L$ is a nilpotent symplectic alternating 
algebra of dimension $4$. We know that $\mbox{dim\, }L/L^{2}\geq 2$. 
If $\mbox{dim\,}L^{2}=2$ then by the reasoning
above, we would have that $\mbox{dim\,}L^{3}=1$ that contradicts Proposition
2.9. By that proposition we neither can have that $\mbox{dim\, }L^{2}=1$. Thus
we must have $L^{2}=0$ and $L$ is abelian.  
\begin{theo} Let $L$ be a nilpotent symplectic alternating algebra of dimension $2n\geq 2$. There
exists an ascending chain of isotropic ideals 
     $$\{0\}=I_{0}<I_{1}<\cdots <I_{n-1}<I_{n}$$
such that $\mbox{dim\,}I_{r}=r$ for $r=0,\ldots ,n$. Furthermore, for $2n\geq 6$, $I_{n-1}^{\perp}$ is
abelian and the ascending chain 

   $$\{0\}<I_{2}<I_{3}<\ldots <I_{n-1}<I_{n-1}^{\perp}<I_{n-2}^{\perp}<\cdots <I_{2}^{\perp}<L$$
is a central chain. In particular $L$ is nilpotent of class at most $2n-3$.
\end{theo}
{\bf Proof}\ \ Starting with the ideal $I_{0}=\{0\}$, we can apply Lemma 
2.3 iteratively to get the required chain 
      $$\{0\}=I_{0}<I_{1}<\ldots <I_{n}.$$
By Lemma 2.7 we have that $I_{2}\leq Z(L)$. By this and Lemma 2.3 we thus have that the chain 
     $$I_{0}< I_{2}< I_{3}< \ldots < I_{n-1}$$
is central in $L$. By Lemma 2.8 it follows that the chain
    $$I_{n-1}^{\perp}< I_{n-2}^{\perp}< \ldots < I_{2}^{\perp}
< I_{0}^{\perp}$$
is also central. It only remains
to see that $I_{n-1}<I_{n-1}^{\perp}$ is central and that $I_{n-1}^{\perp}$ 
is abelian. As $I_{n-1}^{\perp}=I_{n-1}+Fx+Fy$ for some $x,y\in L$ where 
$(x,y)=1$ and as $I_{n-1}$
is isotropic, this follows from Lemma 2.6. $\Box$ \\ \\
{\bf Remark}. When $\mbox{dim\,}Z(L)=r<n$, we can choose our chain such that $I_{r}=Z(L)$. We then get a central chain
      $$I_{0}<I_{r}<I_{r+1}<\cdots <I_{n-1}<I_{n-1}^{\perp}<I_{n-2}^{\perp}<\cdots <I_{r}^{\perp}<L.$$
In particular the class is then at most $2n-3-2(r-2)=2n-2r+1$. \\ \\
{\bf Presentations of nilpotent symplectic alternating algebras}. 
Last proposition tells us a great deal about the structure of nilpotent symplectic alternating algebras. A moments reflection should convince the reader that we can pick a standard
basis $x_{1},y_{1},x_{2},y_{2}, \ldots ,x_{n},y_{n}$ such that
             $$I_{1}=Fx_{n},\ I_{2}=Fx_{n}+Fx_{n-1},\ \cdots ,I_{n}=Fx_{n}+\cdots +Fx_{1},$$
             $$I_{n-1}^{\perp}=I_{n}+Fy_{1},\ I_{n-2}^{\perp}=I_{n}+Fy_{1}+Fy_{2},\ \cdots,\ I_{0}^{\perp}=L=I_{n}+Fy_{1}+\cdots +Fy_{n}.$$
Now let $u,v,w$ be three of the basis elements. Since $I_{n}$ is abelian we have that $(uv,w)=0$ whenever two of these three elements are from
$\{x_{1},\ldots ,x_{n}\}$. The fact that 
               $$\{0\}<I_{1}<\ldots <I_{n}$$
is central also implies that $(x_{i}y_{j},y_{k})=0$ if $i\geq k$. So we only 
need to consider the possible non-zero triples $(x_{i}y_{j},y_{k}),\ \ (y_{i}y_{j},y_{k})$ for $1\leq i<j<k\leq n$. For each triple $(i,j,k)$ with 
$1\leq i<j<k\leq n$, let $\alpha(i,j,k)$ and $\beta(i,j,k)$ be some elements
in the field $F$. We refer to the data 
   $${\mathcal P}:\ \ (x_{i}y_{j},y_{k})=\alpha(i,j,k),\ \ (y_{i}y_{j},y_{k})=
\beta(i,j,k),\ \ \ \ \ 1\leq i<j<k\leq n$$
as a {\it nilpotent presentation}. We have
just seen that every nilpotent symplectic alternating algebra has a 
presentation of this type.  Conversely, given any nilpotent presentation, let 
          $$I_{r}=Fx_{n}+Fx_{n-1}+\cdots +Fx_{n+1-r}$$
and we get an ascending central chain of isotropic ideals $\{0\}=I_{0}<I_{1}<
\ldots <I_{n}$ such that $\mbox{dim\,}I_{j}=j$ for $j=1,\ldots ,n$. By Lemma 2.8
we then get a central chain
   $$\{0\}=I_{0}<I_{1}<\ldots <I_{n}<I_{n-1}^{\perp}<I_{n-2}^{\perp}
  <\ldots <I_{0}^{\perp}=L$$
and thus $L$ is nilpotent. Thus every nilpotent presentation describes a nilpotent
symplectic alternating algebra. \\ \\ 
{\bf Remark}. Notice that there are $2{n\choose 3}$ parameters for these presentations. If $F$ is a finite field this thus gives the value $|F|^{2{n\choose 3}}$
as the upper bound for the number of $2n$-dimensional nilpotent symplectic 
alternating algebras over the field $F$. Armed with this information it is not
difficult to get some good information about the growth of nilpotent symplectic
alternating algebras over a finite field $F$. Let $V$ be a $2n$-dimensional
vector space over $F$ and consider $(\wedge^{3} V)^{*}$. After fixing a standard basis for $V$, each presentation 
of a symplectic alternating algebra corresponds to an element in 
$(\wedge^{3} V)^{*}$. Now let ${\mathcal N}$ be the subset of $(\wedge^{3} V)^{*}$
corresponding to all nilpotent presentations. The number of nilpotent
symplectic alternating algebras of dimension $2n$ is the same as the number
of $Sp(V)$-orbits of $(\wedge^{3} V)^{*}$ consisting of presentations that 
give nilpotent algebras. Suppose these are $u_{i}^{\mbox{Sp}(V)},\ i=1,\ldots ,m$. Then
      $${\mathcal N}\subseteq \cup_{i=1}^{m}u_{i}^{\mbox{Sp}(V)}$$
and thus $|F|^{2{n\choose 3}}=|{\mathcal N}|\leq m\cdot |\mbox{Sp}(V)|
\leq m\cdot |F|^{2n+1\choose 2}$. These calculations show that the number of nilpotent symplectic alternating algebras is
     $$m=|F|^{n^{3}/3+O(n^{2})}.$$
From Proposition 2.9 we know that no term of the lower central series
of a symplectic alternating algebra can be $1$-dimensional. Next proposition
shows that some of terms of the lower central series cannot be 
$2$-dimensional. \
\begin{prop} Let $L$ be a symplectic alternating algebra  we have that
$\mbox{dim\,}L^{m}\not =2$ for $2\leq m\leq 4$. Equivalently 
$Z_{m}(L)$ is not of co-dimension $2$ if $1\leq m\leq 3$.  
\end{prop}
{\bf Proof}. We first prove that $\mbox{dim\,}L^{2}\not
=2$. We argue by contradiction and suppose $\mbox{dim\,}L^{2}=2$. Then 
      $$2=\mbox{dim\,}L^{2}=\mbox{dim\,}Z(L)^{\perp}=
\mbox{dim\,}L-\mbox{dim\,}Z(L).$$
Suppose $L=Z(L)+Fu+Fv$. Then $L^{2}=Fuv$, which contradicts 
$\mbox{dim\,}L^{2}=2$. \\ \\
Next we turn to showing that $\mbox{dim\,}L^{3}\not =2$. We argue by 
contradiction and let $L$ be a counter example of smallest dimension. We first
notice that $Z(L)$ must be isotropic as otherwise $L=I\oplus I^{\perp}$ for
some $2$-dimensional ideal $I=Fu+Fv\leq Z(L)$ where $(u,v)=1$. But then $M=I^{\perp}$ is
a symplectic alternating algebra of smaller dimension where 
$M^{3}=L^{3}$ is of dimension $2$. This however contradicts the minimality of
$L$. We can thus assume that $Z(L)$ is isotropic. Notice that
    $$2=\mbox{dim\,}L^{3}=\mbox{dim\,}Z_{2}(L)^{\perp}=
     \mbox{dim\,}L-\mbox{dim\,}Z_{2}(L).$$
Say, $L=Z_{2}(L)+Fx+Fy$. Then $L^{2}=Z(L)+Fxy$ and, as $Z(L)$ is isotropic
and $xy\in L^{2}=Z(L)^{\perp}$, $L^{2}$ is isotropic. Thus $L^{2}\leq (L^{2})^{\perp}=Z(L)$ and we get the contradiction that $L^{3}=\{0\}$. \\ \\
It now only remains to deal with $L^{4}$. For a contradiction, suppose that
$\mbox{dim\,}L^{4}=2$. Then
    $$2=\mbox{dim\,}L^{4}=\mbox{dim\,}Z_{3}(L)^{\perp}=\mbox{dim\,}L-
\mbox{dim\,}Z_{3}(L).$$
Say $L=Z_{3}(L)+Fu+Fv$. Then $L^{2}\leq Z_{2}(L)+Fuv$ and using the fact
that $Z_{2}(L)\cdot L^{2}=\{0\}$ we get 
     $$L^{2}\cdot L^{2}\leq (Z_{2}(L)+Fuv)L^{2}= Fuv\cdot L^{2}\leq
Fuv\cdot (Z_{2}(L)+Fuv)= F(uv)(uv)=0.$$
Thus $0=(L,L^{2}\cdot L^{2})\Rightarrow (L^{3},L^{2})=0\Rightarrow 
(L^{4},L)=0$, that gives us the contradiction that $L^{4}=\{0\}$. $\Box$ \\ \\
{\bf Example}. Let $L$ be the nilpotent alternating algebra with presentation
(we only list the triples that have non-zero value)
    $$(x_{2}y_{3},y_{4})=1,\ (x_{1}y_{2},y_{3})=1,\ (y_{1}y_{2},y_{4})=1.$$
Then inspection shows that $\mbox{dim\,}L^{5}=2$. The bound $4$ in last 
proposition is therefore the best one. 
\section{The structure of nilpotent symplectic alternating algebras of maximal class}
We have seen previously that nilpotent symplectic alternating algebras of dimension $2n$ have class at most $2n-3$. For every 
algebra of dimension $2n\geq 8$ this bound is attained. As well as 
demonstrating this we will see that the structure of these
algebras of maximal class is very restricted. \\ \\
Let $L$ be a nilpotent symplectic alternating algebra of dimension $2n\geq 8$ 
with an ascending chain of isotropic ideals
    $$\{0\}=I_{0}<I_{1}<\cdots <I_{n},$$ 
where $\mbox{dim\,}I_{j}=j$ for $j=1,\ldots ,n$.  
\begin{theo} Suppose $L$ is of maximal
class. Then
    $$I_{2}=Z_{1}(L),\ I_{3}=Z_{2}(L),\ \ldots ,I_{n-1}=Z_{n-2}(L),$$
   $$I_{n-1}^{\perp}=Z_{n-1}(L),\ I_{n-2}^{\perp}=Z_{n}(L),\ \cdots ,
I_{2}^{\perp}=Z_{2n-4}(L).$$
Furthermore $Z_{0}(L), Z_{1}(L), \ldots ,Z_{2n-3}(L)$ are the unique ideals
of $L$ of dimensions $0,2,3,\ldots ,n-1,n+1,n+2,\ldots ,2n-2, 2n$. 
\end{theo}
{\bf Proof}\ \ Let $J_{0}=\{0\}, J_{1}=I_{2},\ldots ,J_{n-2}=I_{n-1},
J_{n-1}=I_{n-1}^{\perp},J_{n}=I_{n-2}^{\perp},\ldots ,J_{2n-4}=I_{2}^{\perp},
J_{2n-3}=L$.  By Theorem 2.10, the chain $J_{0}<J_{1}<\ldots <J_{2n-3}$ is central.
We argue by contradiction and let $i$ be the smallest
integer between $1$ and $2n-4$ where $J_{i}<Z_{i}(L)$. Let 
$u\in Z_{i}(L)\setminus J_{i}$ and let $k$ be
the smallest integer between $i$ and $2n-4$ such that $u\in J_{k+1}$.
Then
               $$J_{k}<J_{k}+Fu\leq J_{k+1}.$$
If $J_{k+1}/J_{k}$ has dimension $1$ it follows that $J_{k+1}\leq Z_{k}(L)$
and we get the contradiction that the class is at most $2n-4$. We can thus
suppose that $J_{k+1}/J_{k}$ has dimension $2$ and there are two cases to 
consider, either $k=n-2$ or $k=2n-4$. In the former case we have
      $$I_{n-1}<I_{n-1}+Fu\leq I_{n-1}^{\perp}$$
which implies that $I=I_{n-1}+Fu$ is an isotropic ideal of maximal 
dimension $n$. As $u\in Z_{n-2}(L)$, we have that $I_{n-2}<I$ is centralised by $L$. By Lemma 2.8 it follows that $I<I_{n-2}^{\perp}$ is also centralised
by $L$ and we
we get a central series
   $$\{0\}=I_{0}<I_{2}<I_{3}<\ldots <I_{n-2}<I<I_{n-2}^{\perp}<I_{n-3}^{\perp}<
\cdots <I_{2}^{\perp}<I_{0}^{\perp}=L$$
of length $2n-4$ and we get again the contradiction that the class is
less than $2n-3$. Finally
suppose that $k=2n-4$. So we have
      $$I_{2}^{\perp}<I_{2}^{\perp}+Fu<L$$
and $u\in Z_{2n-4}(L)$. Now let $v\in L\setminus (I_{2}^{\perp}+Fu)$. Then 
$L=I_{2}^{\perp}+Fu+Fv$ and $L^{2}=(I_{2}^{\perp}+Fu)L\leq Z_{2n-5}(L)$. Hence $L\leq
Z_{2n-4}(L)$ that again contradicts the assumption that $L$ is of class
$2n-3$. \\ \\
We now want to show that these terms of the upper central series are the unique ideals of dimensions $0,2,3,\ldots ,n-1,n+1,n+2,\ldots ,2n-2,2n$. 
First let $I$ be an ideal
of dimension $2$. By Lemma 2.7 we have that $I\leq Z(L)$ and as we have seen
that $Z(L)$ has dimension $2$, it follows that $I=Z(L)$. Now suppose
that for some $2\leq k\leq n-2$ we know that $Z_{k-1}(L)$ is the only
ideal of dimension $k$. Let $I$ be an
ideal of dimension $k+1$. As $L$ is nilpotent we have that $I$ contains
a ideal $J$ of dimension $k$. By the induction hypothesis we have that 
$J=Z_{k-1}(L)$ and as $I/J$ is of dimension $1$ we have that $I\leq Z_{k}(L)$.
We have that $Z_{k}(L)$ has dimension $k+1$ and thus 
$I=Z_{k}(L)$. We have thus seen that there are unique ideals of dimensions
$0,2,3,\ldots ,n-1$. Now let $I$ be an ideal of dimension $i\in\{
n+1,n+2,\ldots ,2n-2,2n\}$. Then $I^{\perp}$ is an ideal whose dimension is
in $\{0,2,3,\ldots ,n-1\}$. By what we have just seen $I^{\perp}$ is unique
and thus $I$ as well. $\Box$ \\ \\
{\bf Remarks} (1) In particular it follows that $Z_{k}(L)^{\perp}=
Z_{2n-3-k}(L)$ for $0\leq k\leq 2n-3$. \\ \\
(2) As $L^{k}=Z_{k-1}(L)^{\perp}$, it follows that $L,L^{2},\ldots ,L^{2n-2}$
are the unique ideals of dimensions $2n,2n-2,2n-3,\ldots , n+1,n-1,n-2,\ldots
,2,0$. Also 
           $$L^{k}=Z_{k-1}(L)^{\perp}=Z_{2n-k-2}(L).$$
{\bf Remark}. Let $L$ be any nilpotent symplectic alternating algebra of dimension $2n\geq 6$ with the property that $\mbox{dim\,}Z(L)=2$. Notice that $Z(L)$ must be isotropic since otherwise we would have a $2$-dimensional symplectic subalgebra $I$ within $Z(L)$ and we would get a direct sum $I\oplus I^{\perp}$ of two symplectic alternating algebras. As $I^{\perp}$ has non-trivial center this would contradict the assumption that $Z(L)$ is $2$-dimensional. Now $L$ has rank $2$. Suppose it is generated by $x,y$. Then $L$ is generated by $x,y,xy$ modulo $L^{3}$ and thus $\mbox{dim\,}Z_{2}(L)=\mbox{dim\,}(L^{3})^{\perp}=\mbox{dim\,}L-\mbox{dim\,}L^{3}=3$. \\ \\
{\bf The complete list of ideals of $L$}. We have seen that there is a
unique ideal of dimension $k$ for any 
$0\leq k\leq 2n$ apart from $k=1, k=n$ and $k=2n-1$. Let us now turn to
the remaining dimensions. Now every ideal of dimension
$1$ is contained in $Z(L)$ and conversely every subspace of dimension $1$
in $Z(L)$ is an ideal. \\ \\
Next consider an ideal $I$ of dimension $2n-1$. Then $I^{\perp}$ is an
ideal of dimension $1$ and is thus any subspace of dimension $1$ such that
    $$\{0\}<I^{\perp}<Z(L)$$
Equivalently, $I$ is any subspace of dimension $2n-1$ such that
    $$L^{2}=Z(L)^{\perp}<I<\{0\}^{\perp}=L.$$
Finally consider an ideal $I$ of dimension $n$. Since $L$ is nilpotent there exists
an ideal $J$ of dimension $n+1$ containing $I$. By last theorem we
have that $J=L^{n-1}=Z_{n-2}(L)^{\perp}$. Also $I$ contains an ideal of dimension
$n-1$ that we know is $Z_{n-2}(L)$. Thus 
      $$Z_{n-2}(L)<I<Z_{n-2}(L)^{\perp}.$$
We also know from our previous work that $Z_{n-2}(L)$ is an isotropic ideal
of dimension $n-1$. $I$ is thus an isotropic ideal of the form 
            $$Z_{n-2}(L)+Fu$$
For some $u\in Z_{n-2}(L)^{\perp}\setminus Z_{n-2}(L)$. Conversely, as 
$Z_{n-2}(L)^{\perp}L\leq Z_{n-2}(L)$ we have that for any intermediate subspace $I$
of dimension $n$ between $Z_{n-2}(L)$ and $Z_{n-2}(L)^{\perp}$, $I$ is an ideal. \\ \\
We thus have a complete picture of the ideals of $L$. \\ \\
We now focus on the characteristic ideals. It turns out that there are as 
well always characteristic ideals of dimension $1,n$ and $2n-1$ when 
$2n\geq 10$. \\ \\
{\bf Remark}. Notice that if $I$ is a characteristic ideal then the ideal
$I^{\perp}$ is also characteristic. To see this let $\phi$ be any automorphism
of the symplectic alternating algebra $L$ 
and let $a\in I^{\perp}$. As $\phi$ is an automorphism we have that 
$\phi(a)\in \phi(I)^{\perp}=I^{\perp}$.  
\begin{theo} Let $L$ be a nilpotent symplectic alternating algebra of dimension
$2n\geq 10$ that is of maximal class. $L$ has a chain of characteristic ideals
   $$\{0\}=I_{0}<I_{1}<\cdots <I_{n}<I_{n-1}^{\perp}<\cdots <
    I_{1}^{\perp}<I_{0}^{\perp}=L$$
where for $0\leq k\leq n$, $I_{k}$ is isotropic of dimension $k$. 
\end{theo}
{\bf Proof} By Theorems 2.10 and 3.1, we know that we can get such a 
chain of ideals where all the ideals apart from $I_{1},I_{n}$ and
$I_{2n-1}$ are characteristic. We want to show that we can choose our chain
such that $I_{1},I_{n}$ and $I_{2n-1}$ are also characteristic. Let
$x_{1},y_{1},\ldots ,x_{n},y_{n}$ be a standard basis such that 
         $$I_{k}=Fx_{n}+Fx_{n-1}+\cdots +Fx_{n+1-k}$$
for $1\leq k\leq n$. Then $I_{4}I_{2}^{\perp}=Fx_{n-3}y_{n-2}$ is a characteristic
ideal. We claim that this is non-trivial. Otherwise
$x_{n-3}y_{n-2}=0$ and then $(x_{n-3}u,y_{n-2})=0$ for all $u\in L$ that implies
that $x_{n-3}L\leq Fx_{n}+Fx_{n-1}$ and we get the contradiction that
$x_{n-3}\in Z_{2}(L)=I_{3}$. Thus we have got a characteristic ideal 
of dimension $1$, namely $I_{4}I_{2}^{\perp}=Z_{3}(L)\cdot L^{2}$. Notice that
we are assuming here that $n\geq 5$. From this
we get that $(I_{4}I_{2}^{\perp})^{\perp}$ is a characteristic ideal of
dimension $2n-1$. \\ \\
It remains to find a characteristic ideal of dimension $n$. We know that
$L^{n}=Fx_{n}+Fx_{n-1}+\cdots +Fx_{2}$, 
$L^{n-1}=I_{n-1}^{\perp}=Fx_{n}+\cdots +Fx_{1}+Fy_{1}$ and $L^{n-2}=I_{n-2}^{\perp}=
I_{n-1}^{\perp}+Fy_{2}$. As $L^{n-1}=L^{n-2}\cdot L$ it follows that 
$L^{n}+Fx_{1}+Fy_{1}=(L^{n-1}+Fy_{2})L$ and thus
       $$L^{n}+Fx_{1}+Fy_{1}=L^{n}+y_{2}L.$$
Thus there exist $u,v\in L$ such that $y_{2}u+L^{n}=x_{1}+L^{n}$
and $y_{2}v+L^{n}=y_{1}+L^{n}$. Then 
  $$(y_{2}u,x_{1})=0, (y_{2}u,y_{1})\not =0, (y_{2}v,x_{1})\not =0,
   (y_{2}v,y_{1})=0.$$
Equivalently 
  $$(x_{1}y_{2},u)=0, (x_{1}y_{2},v)\not =0, (y_{1}y_{2},u)\not =0,
   (y_{1}y_{2},v)=0$$
and this implies that $x_{1}y_{2}, y_{1}y_{2}$ are linearly independent 
(something that will also be useful later). Consider next the $2$-dimensional
characteristic subspace 
    $$L^{n-1}L^{n-2}=Fx_{1}y_{2}+Fy_{1}y_{2}.$$
Notice that $L^{n-1}L^{n-2}\leq I_{n-2}$.  Let $k$ be the smallest positive integer between $1$ and $n-3$ such that $L^{n-1}L^{n-2}\leq I_{k+1}$. Let  
$J=L^{n-1}L^{n-2}\cap I_{k}$. Then $\mbox{dim\,}J=1$ and there is a unique 
one-dimensional
subspace $Fu$ of $Fx_{1}+Fy_{1}$ such that $FuL^{n-2}=J$. Now 
$I=I_{n-1}+Fu$ is the characteristic ideal of dimension $n$ that we wanted. 
Notice that $I=\{x\in I_{n+1}:xI_{n+2}\subseteq J\}$. 
$\Box$ \\ \\
{\bf Remark}. If $L$ is a nilpotent symplectic alternating algebra of dimension
$8$ that is of maximal class then there is no characteristic ideal
of dimension $1$. The reader can convince himself of this by looking at the
classification of these algebras given in the last section. 
\begin{coro} Let $L$ be a nilpotent symplectic alternating algebra of maximal 
class and dimension $2n\geq 10$. The automorphism group of $L$ is nilpotent-by-abelian. 
\end{coro}
{\bf Proof} Consider a chain of characteristic ideals as given in the
last theorem
     $$\{0\}=I_{0}<I_{1}<\ldots <I_{n}<I_{n-1}^{\perp}<I_{n-2}^{\perp}<
\ldots <I_{0}^{\perp}=L.$$
Consider the ordered basis $(x_{n},x_{n-1},\ldots ,x_{1},y_{1},
\ldots, y_{n})$ associated with this chain, that is $I_{k}=Fx_{n}+ Fx_{n-1}+
\cdots +Fx_{n+1-k}$. As the ideals in the chain are all characteristic we
see that the matrix of any automorphism with respect to that ordered basis will be upper triangular. The result
follows. $\Box$ \\ \\
We next move on to  presentations of nilpotent symplectic alternating algebras of 
maximal 
class. Suppose $L$ is any nilpotent symplectic alternating algebra with a presentation
  $${\mathcal P}:\ \ (x_{i}y_{j},y_{k})=\alpha_{ijk},\ \ 
          (y_{i}y_{j},y_{k})=\beta_{ijk}\ \ \ 1\leq i<j<k\leq n .$$
We would like to read from the presentation whether the algebra is of maximal
class. This turns out to be possible.
\begin{theo} Let $L$ be a nilpotent symplectic alternating algebra of dimension
$2n\geq 8$ given by 
some nilpotent presentation ${\mathcal P}$. The algebra is of maximal class if and
only if $x_{n-2}y_{n-1}$, $x_{n-3}y_{n-2}$, $\ldots$ , $x_{2}y_{3}$ are non-zero
and $x_{1}y_{2},y_{1}y_{2}$ are linearly independent. 
\end{theo}
{\bf Proof} Let us first see that these conditions are necessary. Suppose
that $L$ is of maximal class.  In the proof 
of Theorem 3.2 we have already seen that $x_{1}y_{2}$ and $y_{1}y_{2}$ must
be linearly independent. As before we let $I_{k}=Fx_{n}+\ldots +
Fx_{n+1-k}$. As $x_{n-2}\not\in Z(L)$, we have $(x_{n-2}y_{n-1},y_{n})\not
=0$ and thus $x_{n-2}y_{n-1}\not =0$. As the terms of the central chain
       $$I_{0}<I_{2}<I_{3}<\ldots <I_{n-1}<I_{n-1}^{\perp}<I_{n-2}^{\perp}<
\ldots <I_{2}^{\perp}<I_{0}^{\perp}$$
are the terms of the lower central series, we know that $I_{k+1}L=I_{k}$
for $2\leq k\leq n-2$. Thus we have for $3\leq k\leq n-2$ that 
       $$I_{k-1}+Fx_{n-k+1}=(I_{k}+Fx_{n-k})L.$$
From this it follows $I_{k-1}+Fx_{n-k+1}=I_{k-1}+x_{n-k}L$. In particular there 
exists $u\in L$ such that $I_{k-1}+x_{n-k+1}=I_{k-1}+x_{n-k}u$. It follows 
that $0\not =(x_{n-k}u,y_{n-k+1})=-(x_{n-k}y_{n-k+1},u)$. Hence $x_{n-k}y_{n-k+1}$
is non-zero for $3\leq k\leq n-2$. \\ \\ 
Let us then see that these conditions are sufficient. We do this by showing that
$I_{2}=I_{3}L$, $I_{3}=I_{4}L$, $\ldots $, $I_{n-2}=I_{n-1}L$, 
$I_{n-1}=I_{n-1}^{\perp}L$, $I_{n-1}^{\perp}=I_{n-2}^{\perp}L$, $\ldots$,
$I_{3}^{\perp}=I_{2}^{\perp}L$, $I_{2}^{\perp}=I_{0}^{\perp}L$. This is sufficient as
this would imply that $L^{2n-3}=I_{2}\not =0$ and thus $L$ is nilpotent of
class $2n-3$. Firstly as $x_{n-2}y_{n-1}\not =0$ we have that 
$(x_{n-2}y_{n-1},y_{n})\not =0$ and thus $I_{3}L=x_{n-2}L=Fx_{n}+Fx_{n-1}=
I_{2}$. Now suppose that we have already established that $I_{k}=
I_{k+1}L$ for all $2\leq k\leq m$ where $2\leq m\leq n-3$. Then 
       $$I_{m+2}L=(I_{m+1}+Fx_{n-m-1})L=I_{m}+x_{n-m-1}.L$$
As $x_{n-m-1}y_{n-m}\not =0$ we have $(x_{n-m-1}u,y_{n-m})=-
(x_{n-m-1}y_{n-m},u)\not =0$ for some
$u\in L$ and thus $I_{m+2}L=I_{m}+x_{n-m-1}L=I_{m}+Fx_{n-m}=I_{m+1}$. We have thus
established by induction that 
     $$I_{2}=I_{3}L, \ldots ,I_{n-2}=I_{n-1}L.$$
We next show that $I_{n-1}^{\perp}L=I_{n-1}$. As $x_{1}y_{2}\not =0$ we have
that there exist $u\in L$ such that $0\not =(x_{1}y_{2},u)=-(x_{1}u,y_{2})$
and
    $$I_{n-1}^{\perp}L=(I_{n-1}+Fx_{1}+Fy_{1})L=
      I_{n-2}+Fx_{2}=I_{n-1}.$$
Next we show that $I_{n-2}^{\perp}L=I_{n-1}^{\perp}$.     
As $x_{1}y_{2},y_{1}y_{2}$ are
linearly independent there exist $u,v\in L$ such that 
     $$(x_{1}y_{2},u)=0, (x_{1}y_{2},v)\not =0, (y_{1}y_{2},u)\not =0,
   (y_{1}y_{2},v)=0$$
and thus 
    $$(y_{2}u,x_{1})=0, (y_{2}u,y_{1})\not =0, (y_{2}v,x_{1})\not =0,
   (y_{2}v,y_{1})=0.$$
Hence
     $$I_{n-2}^{\perp}L=(I_{n-1}^{\perp}+Fy_{2})L=
            I_{n-1}+Fy_{2}L=I_{n-1}+Fx_{1}+Fy_{1}=I_{n-1}^{\perp}.$$
Now suppose that we have established that $I_{k-1}^{\perp}L=I_{k}^{\perp}$ for 
$m+1\leq k\leq n-1$ where $3\leq m\leq n-2$. As $x_{n-m}y_{n-m+1}\not =0$
it follows that there exists $u\in L$ such that $0\not =(x_{n-m}y_{n-m+1},u)=
(y_{n-m+1}u,x_{n-m})$. Thus  
  $$I_{m-1}^{\perp}L=(I_{m}^{\perp}+y_{n-m+1})L=I_{m+1}^{\perp}+y_{n-m+1}L
     =I_{m+1}^{\perp}+Fy_{n-m}=I_{m}^{\perp}.$$
It now only remains to see that $I_{0}^{\perp}L=I_{2}^{\perp}$. But this follows from $x_{n-2}y_{n-1}\not =0$ that implies that $(y_{n-1}y_{n},x_{n-2})=
(x_{n-2}y_{n-1},y_{n})\not =0$. Thus 
    $$I_{0}^{\perp}L=(I_{2}^{\perp}+Fy_{n-1}+Fy_{n})L=I_{3}^{\perp}+
(Fy_{n-1}+Fy_{n})L=I_{3}^{\perp}+Fy_{n-2}=I_{2}^{\perp}.$$
This finishes the proof. $\Box$ \\ \\
{\bf Remark}. In particular it follows that for each $2n\geq 8$ there exist
a nilpotent symplectic alternating algebra of maximal class. One just needs
to choose the presentation such that the conditions from Theorem 3.4 hold. One 
possibility is
    $$\begin{array}{l}
      {\mathcal P}:\  (x_{n-2}y_{n-1},y_{n})=-1,\ (x_{n-3}y_{n-2},y_{n})=-1,\ 
\cdots ,
   (x_{1}y_{2},y_{n})=-1, \\
\mbox{\ \ \ \ \ \ } (y_{1}y_{2},y_{n-1})=-1.
\end{array}$$     
In fact the conditions are not a strong constraint. In particular the 
values of $(x_{i}y_{j},y_{k}), (y_{i}y_{j},y_{k})$ where $j-i\geq 2$ can
be chosen freely. The number of such triples is $2{n-1 \choose 3}$ that is a 
polynomial in $n$ of degree $3$ with leading coefficient $1/3$. Let $F$ be
any finite field. By a similar
argument as we used for determining the growth of nilpotent symplectic 
alternating algebras 
we see that the number $m(n)$ of nilpotent symplectic alternating algebras
of maximal class satisfies 
      $$m(n)=|F|^{n^{3}/3+O(n^{2})}.$$
{\bf Remark}. (1) Let $L$ be a nilpotent symplectic alternating algebra 
of dimension $2n\geq 10$ that is of maximal class and consider a chain $\{0\}=I_{0}<\ldots <I_{n}$ of characteristic
ideals where $I_{k}$ is of dimension $k$. We have, for $4\leq m\leq n-1$,
    $$I_{m}I_{m-2}^{\perp}=Fx_{n+1-m}y_{n+2-m}$$
and thus we get that $Fx_{2}y_{3}, Fx_{3}y_{4}, \ldots ,Fx_{n-3}y_{n-2}$ are 
one-dimensional characteristic subspaces of $L$. Also
             $$I_{n-1}^{\perp}I_{n-2}^{\perp}=Fx_{1}y_{2}+Fy_{1}y_{2}$$
is a characteristic subspace. So is $I_{n}^{\perp}I_{n-2}^{\perp}=Fx_{1}y_{2}$. \\ \\
(2) If $V$ is a characteristic subspace of dimension $d$ then we get 
a chain of characteristic subspaces 
   $$V\cap I_{1}\subseteq V\cap I_{2}\subseteq
\ldots \subseteq V\cap I_{n}\subseteq V\cap I_{n-1}^{\perp}\subseteq \ldots
\subseteq V\cap I_{0}^{\perp}=V.$$
Thus there is a chain of characteristic subspaces $V_{1}<V_{2}<
\ldots <V_{d}$ where $v_{i}$ is of dimension $i$. 
\section{Nilpotent algebras of dimension $2n\leq 8$}
The classification of the nilpotent symplectic alternating algebras of dimensions at most $8$ is implicit in [3] although this is not done explicitly and 
the context there is a more general setting. To demonstrate the machinery that we have developed we will offer a much shorter approach here. The classification of algebras of dimension $10$ is far more challenging and will be dealt with 
in a sequel to this paper. Through this section we will be working with an 
arbitrary field $F$. \\ \\
We have observed earlier that nilpotent symplectic alternating algebras
of dimensions $2$ or $4$ must be abelian. 
\subsection{Algebras of dimension 6}
Let $L$ be a non-abelian nilpotent symplectic alternating algebra of dimension $6$ with a nilpotent presentation ${\mathcal P}$. There are at most two non-zero
triple values 
                   $$(x_{1}y_{2},y_{3})=a,\ \ (y_{1}y_{2},y_{3})=b.$$
As $L$ is non-abelian, one of these must be non-zero and, by replacing 
$x_{1},y_{1}$ by 
$-y_{1},x_{1}$ if necessary, we can assume that
$b\not =0$. Replacing then $x_{3},y_{3}$ by $bx_{3},\frac{1}{b}y_{3}$
implies that we can further assume that $(y_{1}y_{2},y_{3})=1$. Finally
replacing $x_{1},y_{1}$ by $x_{1}-ay_{1},y_{1}$ and we can also assume that
$(x_{1}y_{2},y_{3})=0$. Apart from the abelian algebra, there is thus only one
algebra of dimension $6$ with presentation 
    $${\mathcal P}_{1}:\ \ (y_{1}y_{2},y_{3})=1.$$
(We will normally only write down those triples where the value is non-zero).
\subsection{Algebras of dimension 8}
First suppose that $Z(L)$ is not isotropic. We can then choose our standard
basis such that $I=Fx_{4}+Fy_{4}\subseteq Z(L)$ and we get a direct sum
$I\oplus I^{\perp}$ of symplectic alternating algebras of dimensions $2$ and
$6$. From 4.1 we then know that apart from the abelian algebra, there is
only one such algebra $L_{2}=Fx_{4}+Fx_{3}+Fx_{2}+Fx_{1}+Fy_{1}+
Fy_{2}+Fy_{3}+Fy_{4}$ with presentation 
    $${\mathcal P}_{2}:\ \ (y_{1}y_{2},y_{3})=1.$$
We then turn to the situation where $Z(L)$ is isotropic. Let us first see 
that $\mbox{dim\,}Z(L)\not =4$. We argue by contradiction and suppose that
$\mbox{dim\,}Z(L)=4$. Pick a standard basis such that $Z(L)=Fx_{4}+Fx_{3}+
Fx_{2}+Fx_{1}$. Now $L$ is not abelian and thus $(y_{i}y_{j},y_{k})\not =0$
for some $1\leq i<j<k\leq 4$. Without loss of generality, we can suppose
that $(y_{1}y_{2},y_{3})=1$. Suppose now that $(y_{1}y_{2},y_{4})=a$,
$(y_{2}y_{3},y_{4})=b$ and $(y_{3}y_{1},y_{4})=c$. Let $\bar{y_{4}}=
y_{4}-by_{1}-cy_{2}-ay_{3}$. Inspection shows that $\bar{y_{4}}$ is orthogonal
to $L^{2}=Fy_{1}y_{2}+Fy_{2}y_{3}+Fy_{3}y_{1}+F\bar{y_{4}}y_{1}+
F\bar{y_{4}}y_{2}+F\bar{y_{4}}y_{3}$. Thus $\bar{y_{4}}\in (L^{2})^{\perp}=Z(L)$ and
we get the contradiction that $\mbox{dim\,}Z(L)\geq 5$. Thus we have shown
that $\mbox{dim\,}Z(L)\not =4$ and as $\mbox{dim\,}Z(L)$ is always at least
$2$, we have two cases to consider: $\mbox{dim\,}Z(L)=3$ and $\mbox{dim\,}Z(L)
=2$. \\ \\
\underline{\mbox{Dim\,}$Z(L)=3$}. We can choose the standard basis such that 
$Z(L)=Fx_{4}+Fx_{3}+Fx_{2}$ and $L^{2}=Z(L)^{\perp}=Fx_{4}+Fx_{3}+Fx_{2}+Fx_{1}+
Fy_{1}$. By Theorem 2.10, we know that $L^{3}=L^{2}L\subseteq Z(L)$ and by 
Proposition 2.11 we must then have $L^{3}=Z(L)$. As $x_{1}\not\in Z(L)$, we must
have $(x_{1}y_{i},y_{j})\not =0$ for some $2\leq i<j\leq 4$. Without loss 
of generality $(x_{1}y_{2},y_{3})\not =0$. By replacing $y_{4},y_{1}$
by $y_{4}-ax_{1}, y_{1}+ax_{4}$ for a suitable $a$, we can assume that 
$(y_{2}y_{4},y_{3})=0$. Let $V=Fy_{2}+Fy_{3}+Fy_{4}$. Now $(y_{2}y_{3},y_{4})=0$
and $L^{2}=Z(L)+Fx_{1}+Fy_{1}$. As $L=L^{2}+V$ it follows that $V^{2}=
Fx_{1}+Fy_{1}$ and as $V^{2}$ is not isotropic we must have that some two
of $y_{2}y_{3}, y_{4}y_{3}, y_{2}y_{4}$ are not isotropic. Without loss of 
generality we can suppose that these are $y_{2}y_{3}$ and $y_{4}y_{3}$. 
By replacing $y_{4},x_{4}$ by $ay_{4},\frac{1}{a}x_{4}$ for a suitable 
$a\in F$, we can furthermore assume that $(y_{2}y_{3},y_{4}y_{3})=1$. 
Thus 
              $$Fx_{1}+Fy_{1}=V^{2}=Fy_{2}y_{3}+Fy_{4}y_{3}$$
and $y_{2}y_{4}=ay_{2}y_{3}+by_{4}y_{3}$ for some $a,b\in F$. It follows that
$(y_{2}+by_{3})(y_{4}-ay_{3})=0$. Now replace $y_{2},y_{4},x_{3}$ by 
$y_{2}+by_{3}, y_{4}-ay_{3}, x_{3}-bx_{2}+ax_{4}$ and then replace $x_{1},y_{1}$
by $y_{2}y_{3}, y_{4}y_{3}$. It follows that we get a new standard basis
where 
       $$y_{2}y_{3}=x_{1},\ y_{4}y_{3}=y_{1},\ y_{2}y_{4}=0.$$
This implies that the only non-zero triples are $(y_{1}y_{2},y_{3})=1$
and $(x_{1}y_{3},y_{4})=1$. There is thus only one possible candidate
here, the algebra $L_{3}$ with presentation 
   $${\mathcal P}_{3}:\ \ (y_{1}y_{2},y_{3})=1,\ (x_{1}y_{3},y_{4})=1.$$
Conversely, one sees by inspection that $Z(L_{3})=Fx_{4}+Fx_{3}+Fx_{2}$ and
this candidate is a genuine example with $\mbox{dim\,}Z(L)=3$. \\ \\
\underline{$\mbox{Dim\,}Z(L)=2$}. We know that the class of $L$ is at most
$2\cdot 4-3=5$ and thus $L^{5}\leq Z(L)$. Let $k$ be the smallest positive 
integer $2\leq k\leq 5$ such that $L^{k}\leq Z(L)$. As $\mbox{dim\,}L^{k}\leq 2$, it follows from  Proposition 2.11 that $k=5$. Hence $L$ is of maximal class
and by Theorem 3.1 we can choose our standard basis such that, we get ideals
$I_{k}=Fx_{n}+\cdots +Fx_{n+1-k}$, 
$k=0,\ldots ,n$ where
    $$\{0\}=I_{0}<I_{1}<\ldots <I_{4}=I_{4}^{\perp}<I_{3}^{\perp}<\ldots
<I_{0}^{\perp}=L$$
is a central series with $I_{2}=Z(L)=L^{5}$, $I_{3}=Z_{2}(L)=L^{4}$,
$I_{3}^{\perp}=Z_{3}(L)=L^{3}$ and $I_{2}^{\perp}=Z_{4}(L)=L^{2}$. By Theorem 3.4 we
furthermore have that $x_{1}y_{2},y_{1}y_{2}$ are linearly independent and
thus a basis for $Z(L)=Fx_{4}+Fx_{3}$. We can now pick our standard basis such that $x_{1}y_{2}=x_{4}$ and $y_{1}y_{2}=x_{3}$. As $x_{2}\not\in Z(L)$, we also
have that $(x_{2}y_{3},y_{4})\not =0$. This means that we have the non-zero
triples $(x_{1}y_{2},y_{4})=1$, $(y_{1}y_{2},y_{3})=1$ and $(x_{2}y_{3},y_{4})=
r\not =0$. The only remaining triples that are possibly non-zero are
$(x_{1}y_{3},y_{4}), (y_{1}y_{3},y_{4})$ and $(y_{2}y_{3},y_{4})$. Replacing 
$x_{1}$,$y_{1}$, and $y_{2}$ by $x_{1}-ax_{2}$, $y_{1}-bx_{2}$ and $y_{2}+ay_{1}-bx_{1}-cx_{2}$ for suitable $a,b,c\in F$, we can assume that these extra triples are
zero. We can thus choose our basis so that our algebra $L(r)$ has presentation
$${\mathcal P}(r):\ \ (x_{2}y_{3},y_{4})=r,\ (x_{1}y_{2},y_{4})=1,\ 
    (y_{1}y_{2},y_{3})=1.$$
We finally need to sort out when, for $r,s\in F^{*}=F\setminus \{0\}$, $L(r)$
and $L(s)$ are isomorphic. We will see that this happens if and only if
$r/s\in  (F^{*})^{3}$. To see that this is a sufficient condition, suppose
we have an algebra $L$ that has a presentation ${\mathcal P}(r)$ with 
respect to some standard basis $x_{1},y_{1},x_{2},y_{2},x_{3},y_{3},x_{4},y_{4}$. 
Suppose that $s=a^{3}r$ for some $a\in F^{*}$. Let
$\bar{x_{1}}=x_{1}$, $\bar{y_{1}}=y_{1}$, $\bar{x_{2}}=ax_{2}$,
$\bar{y_{2}}=\frac{1}{a}y_{2}$, $\bar{x_{3}}=\frac{1}{a}x_{3}$,
$\bar{y_{3}}=ay_{3}$, $\bar{x_{4}}=\frac{1}{a}x_{4}$ and $\bar{y_{4}}=ay_{4}$.
Inspection shows that $L$ has presentation ${\mathcal P}(s)$ with respect
to the new basis. Hence $L(s)\cong L(r)$ when 
$r/s\in (F^{*})^{3}$. It remains to see that the condition
is also necessary. Consider the algebra $L(r)$ and take an arbitrary new
standard basis $\bar{x_{1}}$, $\bar{y_{1}}$, $\bar{x_{2}}$,
$\bar{y_{2}}$, $\bar{x_{3}}$, $\bar{y_{3}}$, $\bar{x_{4}}$, $\bar{y_{4}}$ such 
that $L(r)$ satisfies the presentation ${\mathcal P}(s)$ for some $s\in F^{*}$.
We want to show that $s/r\in (F^{*})^{3}$. Now 
     $$\bar{y_{3}}=ay_{3}+by_{4}+u,\ \bar{y_{4}}=cy_{3}+dy_{4}+v$$
for some $u,v\in L^{2}$ and $a,b,c,d\in F$ where $ad-bc\not =0$. As 
$\mbox{dim\,}L^{2}-\mbox{dim\,}L^{3}=1$ it follows readily that $L^{2}L^{2}\leq L^{4}$ and
it follows that 
\begin{eqnarray*}
   \bar{y_{3}}\bar{y_{4}}\bar{y_{3}} & = & (ay_{3}+by_{4})(cy_{3}+dy_{4})
           (ay_{3}+by_{4})+w, \\
  \bar{y_{3}}\bar{y_{4}}\bar{y_{4}} & = & (ay_{3}+by_{4})(cy_{3}+dy_{4})
         (cy_{3}+dy_{4})+z,
\end{eqnarray*}
for some $w,z\in L^{4}$. As $L^{6}=0$ we have that $L^{4}$ is orthogonal
to $L^{3}$ and thus in the following direct calculations we can omit
$w$ and $z$. We have 
  $$-s^{2}=(\bar{y_{3}}\bar{y_{4}}\bar{y_{3}},
         \bar{y_{3}}\bar{y_{4}}\bar{y_{4}})=-(ad-bc)^{3}r^{2}.$$
Hence $s/r\in (F^{*})^{3}$. \\ \\
{\it Acknowledgement}. We thank the referee for a number of insightful remarks that led to an improvement of the exposition as well as few 
simplifications.

\bibliographystyle{plainurl}
\bibliography{biblio}

\end{document}